\documentclass[11pt,a4paper]{article}
\usepackage{amsmath}

\title{Polynomial representation of\\
{Fermat's Last Theorem}
\author{Daniele De Pedis \\
Istituto Nazionale di Fisica Nucleare - Roma1 - Italy\\
email: daniele.depedis@roma1.infn.it}}
\date{}
\begin{document}
\maketitle

\begin{abstract}We propose a new approach at Fermat\rq{}s Last Theorem (FLT) solution: for each FLT equation we associate a polynomial of the same degree. The study of the roots of the polynomial allows us to investigate the FLT validity. This technique, certainly within the reach of Fermat himself, allows us infer that this is the \textit{marvelous proof} that Fermat claimed to have. 
\end{abstract}

\setcounter {section}{0}
\section{Introduction}
In 1637 Pierre de Fermat wrote in the margins of a copy of Diophant’s Arithmetical, the book where he used to write many of his famous theories~[1]:
\begin{quotation}
\textit{"It is impossible to separate a cube into two cubes or a fourth power into two fourth powers or, in general, all the major powers of two as the sum of the same power. I have discovered a truly marvelous proof of this theorem, which can\rq{}t be contained in the too narrow page margin".}
\end{quotation}

In other words, the previous expression can be condensed into:  \\
the equation:                                                                                             
\begin{equation} \label{eq.1} 
{A^n+B^n=C^n} 
\end{equation}
 has no solutions, other than trivial ones\footnote{The trivial solution is a solution with at least one of the integers A, B and C equal to zero}, for any value of $A, B, C, n$ integers and n$>$2.

The equation (\ref{eq.1}) is known as Fermat's Last Theorem (FLT). Last, not because it was the last work of Fermat in chronological sense, but because it has remained for over 350 years the Fermat's theorem never solved. In fact, also the same Fermat, although stating the unsolvability of (1) he never provided a complete demonstration (maybe lost) but has left his proof limited only to the case n = 4. In reality, therefore, it would be more correct to talk about Fermat\rq{}s conjecture.

Today, many mathematicians are of the opinion that Fermat was wrong and that he had not a real full demonstration. Others think that Fermat had such proof, or at least that he had guessed the road, but, as was his custom, he was so listless that such evidence went lost.
In any case, as you wish to take a position, the fact remains that for over 350 years all the greatest mathematicians have tried to find such evidence without success.

Only in 1994, after seven years of complete dedication to the problem, Andrew Wiles, who was fascinated by the theorem that as a child dreamed to solve, finally managed to give a demonstration. Since then, we might refer to (\ref{eq.1}) as Fermat's theorem.

However, Wiles used elements of mathematic and modern algebra [2] that Fermat could not know: the demonstration that Fermat claimed to have, if it were correct, then must be so different. 

In this paper we\rq{}ll try to give our contribution proposing a demonstration of Fermat\rq{}s Last Theorem using a technique certainly within the reach of Fermat himself, and then infer that this is the {\itshape marvelous proof} that Fermat claimed to have. In agreement to the supposed Fermat\rq{}s knowledge, we'll also avoid using procedure and notations proper of modern algebra.

\section{First considerations [3][4]}

Before to go in deep in the proof, we make some well known\footnote{See [3] pag.2}  considerations relating to (1).

a)	According to the usual spoken, to say that Fermat's theorem is true is equivalent to saying that (1) is never verified. Nevertheless, the trivial solution is a true solution that we have to consider as we\rq{}ll see later.

b)	It is sufficient to prove (1) be true for the exponent n = 4 and for every n = odd prime. 
As mentioned the case of n = 4 was proved directly by Fermat.

c)	A, B, C must be such that their Greater Common Divider $(GCD)$ is the unit when taken in pairs, i.e.: \\
$GCD(A,B)=GCD(A,C)=GCD(B,C)=1$ and also $GCD(A,B,C)=1$.

d)	Important corollary to the previous property is that the three variables $A, B, C$  can\rq{}t all have the same parity and, moreover, only one can be even following this scheme:

\small
\begin{center}
\begin{tabular}{|c|c|c|c|}
\hline
& A & B & C \\
\hline
1 & odd & odd & even \\
\hline
2 & odd & even & odd \\
\hline
3 & even & odd & odd \\
\hline
\end{tabular} ~ Tab.1 
\end{center}

%\begin{center}
%~Tab.1 
%\end{center}
\normalsize

e)	Another important corollary of c) is: \\
          $GCD (A+B, C-A) = 1$\\                                                                               
          $GCD (A+B, C-B) = 1$\\
          $GCD (C-A, C-B)  = 1$

\section{Demonstration of FLT}
Here we consider the case 1 of Tab.1, that is $A$ and $B$ both odd. \\
The cases 2 and 3 in Tab.1 will be discussed in Appendix A.\\
Let  
\begin{equation} 
      D = C-A = odd\ integer
\end{equation}
      $$E = C-B = odd\ integer$$
then, by the considerations at previous point e), we have $GCD (D,E) = 1$. \\
From (1) and (2) we obtain 
\begin{equation} \label{eq.3} 
{A^n+B^n=(C-D)^n+(C-E)^n=C^n} 
\end{equation}
where  $n$ = prime number $\geq$ 3, then developing the powers of binomials

\begin{equation}
(C-D)^n=\sum_{k=0}^{n}(-1)^k  {n \choose k}  D^k C^{n-k}
=C^n+\sum_{k=1}^{n}(-1)^k {n \choose k} D^k C^{n-k}
\end{equation}
\begin{equation}
(C-E)^n=\sum_{k=0}^{n}(-1)^k  {n \choose k} E^k C^{n-k}
=C^n+\sum_{k=1}^{n}(-1)^k  {n \choose k} E^k C^{n-k}
\end{equation}
where, for convenience, we have released the first term under the sign of summation. 
Substituting (4) and (5) in (3) we obtain the fundamental relationship
\begin{equation}
P_{(C,n)}=C^n+ \sum_{k=1}^{n}(-1)^k {n \choose k} \left ( D^k +E^k\right) C^{n-k}=0
\end{equation}

The (6) is the expression of a polynomial (that we call \textit{associated polynomial}) in the unknown $C$, complete, of degree n and with integer coefficients. The fundamental theorem of algebra assures us that there are n roots of (6) which can be: separate, (partially) overlapping, integer, irrational or complex conjugates\footnote{The polynomial in (6) being monic and with all integer coefficients cannot have rational no-integer roots [5]. Moreover, as well as in case of complex roots, the irrational roots must appear in conjugate pairs, that is, if $a+\sqrt{b}$ is an irrational root of (6) then also $a-\sqrt{b}$ is a root, where $a$ and $b$ are integer numbers and $\sqrt{b}$ is irrational. See appendix C}.

Whatever the type of roots, what interests us is the existence of possible integer roots of (6), in fact, given any integers $D, E$ and $n$, if we could find at least one integer solution, other than the trivial one, into full set \{$\Gamma_i$\} of its roots then, using the relations (2), we could get back $A$ and $B$ and disprove Fermat\rq{}s Theorem.

In other words, if we can prove that the (6) has no integer solutions in $C$, anyhow chosen $D, E$ and $n$, then we can never disprove the theorem and therefore Fermat was right, that is the (1) has no solution in the ring of integers.
Equivalently we can state the following
\newtheorem{teorema}{Lemma}
\begin{teorema}
given any integers D, E, and n such that $GCD (D, E)=1$ and  n$\geq$3, showing that (6) does not admit any integer solution, other than the trivial one, for the unknown variable C is equivalent to prove that the Fermat\rq{}s Last Theorem is true.
\end{teorema}

The fundamental theorem of algebra assures us that the polynomial in (6) can be expressed as
\begin{equation}
P_{(C,n)}=\prod_{i=1}^{n}\left ( C-\Gamma_i \right)=0
\end{equation}
where $\Gamma_i$  are the roots of (6).

In order that (6) and (7) are equal, it is necessary and sufficient that the coefficients of the terms of same degree in $C$ are equal.
Expanding (6) and (7) in their terms, we get:

\begin{eqnarray*} 
P_{(C,n)} & = & C^n-{n \choose 1} \left (D+E\right) C^{n-1}+{n \choose 2} \left ( D^2 +E^2\right) C^{n-2}-... \nonumber \\
                 & + & {n \choose {n-1}} \left ( D^{n-1} +E^{n-1}\right) C- {n \choose n}\left ( D^n +E^n\right)=0 \hspace{1.6cm} (6a) \\
 \end{eqnarray*}                
\begin{eqnarray*}
P_{(C,n)}& = & \prod_{i=1}^{n}\left ( C-\Gamma_i \right)=\left ( C-\Gamma_1 \right)\left ( C-\Gamma_2 \right) ... \left ( C-\Gamma_{n-1} \right)\left ( C-\Gamma_n\right)=\nonumber \\
                & =&C^n-\left (\sum_{{i_1}=1}^{n}\Gamma_{i_1}\right)C^{n-1}+\left (\sum_{1\leq{i_1}<{i_2}}^{n}\Gamma_{i_{1}}\Gamma_{i_{2}}\right)C^{n-2}-...\nonumber \\
                 & +& \left (\sum_{1\leq{i_1}<{i_2}<...<{i_{n-1}}}^{n}\Gamma_{i_{1}}\Gamma_{i_{2}}...\Gamma_{i_{n-1}}\right)C-\left   (\Gamma_1\Gamma_2...\Gamma_n\right)=0 \hspace{1.2cm} (7a) \\
\end{eqnarray*}

The (7a) shows that the development of (7) leads to an expression which is the sum of terms with decreasing powers in $C$ and whose ${(n-k)th}$ coefficient is related to the sum of all possible combinations, without repetition, of the n roots taken k-at-a-time.

Equating the coefficients of terms of equal degree in (6a) and (7a), we arrive at the following fundamental system of equations (also known as Viete\rq{}s formula):
$$
{\mathcal \ } 
 \left\{
\begin{array}{llll}
%\begin{alignedat}
\displaystyle {n \choose 1}\left (D+E\right)={\sum_{{i_1}=1}^{n}}\Gamma_{i_1}=\\
 \displaystyle \hspace{1.0cm}=\Gamma_1+(\Gamma_2+...+\Gamma_{n-1}+\Gamma_n ) \hspace{2.3cm}=\Gamma_1+t_1\hspace{1.3cm}(8a) \\
\\
\displaystyle {n \choose 2}\left (D^2+E^2\right) ={{\sum_{1\leq{i_1}<{i_2}}^{n}}\left(\Gamma_{i_1}\Gamma_{i_2}\right)=}\\
\\
\hspace{1.0cm}=\Gamma_1(\Gamma_2+\Gamma_3+...+\Gamma_n)+ 
\displaystyle {{\sum_{2\leq{i_2}<{i_3}}^{n}}\left(\Gamma_{i_2}\Gamma_{i_3}\right)}
\hspace{0.2cm}=\Gamma_1t_1+t_2 \hspace{1.0 cm}(8b)\\
\\
\cdots \cdots\\
\\
\displaystyle {n \choose {n-1}}\left (D^{n-1}+E^{n-1}\right) =
\displaystyle
{\sum_{1\leq{i_1}<{i_2}<...<{i_{n-1}}}^{n}}\left(\Gamma_{i_1}\Gamma_{i_2}...\Gamma_{i_{n-1}}\right)= \\
\\
 \displaystyle \hspace{1.0cm}
=\Gamma_1\left( {\sum_{2\leq{i_2}<...<{i_{n-1}}}^{n}}\Gamma_{i_2}\Gamma_{i_3}...\Gamma_{i_{n-1}}\right)+\\
\\
\displaystyle  \hspace{1.5cm}+\Gamma_2\Gamma_3...\Gamma_{n-1}\Gamma_n
\hspace{4.1cm}=\Gamma_1t_{n-2}+t_{n-1}\hspace{0.2cm}(8c)\\
\\
\displaystyle {n \choose n}\left (D^n+E^n\right)=\Gamma_1\left(\Gamma_2...\Gamma_{n-1}\Gamma_n 
\right)\hspace{2.3cm}=\Gamma_1t_{n-1}\hspace{1.4cm}(8d)
%\end{alignedat}
\end{array}
\right.
$$
\\
where, $\displaystyle t_1 = (\Gamma_2+...+ \Gamma_{n-1} + \Gamma_n), \ 
t_2 ={\sum_{2\leq{i_2}<{i_3}}^{n}}\left(\Gamma_{i_2}\Gamma_{i_3}\right)=(\Gamma_2\Gamma_3+\dots+
\Gamma_2\Gamma_n+\dots+\Gamma_{n-1}\Gamma_n)$  and so on, and in particular $\displaystyle t_{n-1} = (\Gamma_2 \Gamma_3...\Gamma_{n-1} \Gamma_n)$. Moreover, let, without loss of generality, $\Gamma_1$ be the integer trivial root, we will show that it is the only possible integer root.
From equations (8) follows the important

\begin{teorema}
If D and E have the same parity then the terms on the right side of each equation in (8) must have an even integer value.
\end{teorema}

We will see that the condition $D$ and $E$ both odd is incompatible to fulfill all the relations (8).

According to Lemma 1, to prove the FLT, we have to show that the associate polynomial admits one, and only one, integer root and it is the trivial solution\footnote{Here, we impose at the trivial solution only the constraint to be integer. It is straightforward to verify into (6) that cases in which ABC=0 imply $P_{(C,n)}=0$.}.\\
\\
{\bfseries Proof:}                                            
We begin observing that in the equation (6a) all terms, except the first one, contain the even factor $(D^i+E^i)$, therefore, if some integer root exists then it must have an even value.
\\

Now, by Lemma 2, the term on the right side of each equation in the system (8) must be an even integer value, then: 
\begin{itemize}
\item The case in which any of the $t_i $ is non-integer is obviously ruled out. 
\item  $\Gamma_1$ being an even root of (6a) then also $t_1$ and all the $t_i$ must be even integers\footnote{The recursive form of the (8) implies the propagation of the $t_i\rq{}s$ parity along all the equations. In fact in (8a) $\Gamma_1$ and $t_1$, in order their sum is even, must have the same parity then, in (8b) also $t_2$ must have the $t_1$  parity, and so on. On the other hand, $\Gamma_1$ can not be odd, because otherwise the term $\Gamma_1t_{n-1}$ in (8d) would also be odd in contradiction with Lemma 2.}.
\item The left side of equation (8c)\footnote{The (8c) in general will be the penultimate equation of any system with n=p equations. On the left side of this equation there is always the sum of two even powers (i.e. n-1) of odd terms. On the right side, due to the construction procedure  of Viete\rq{}s formulas, there will be always the sum of terms made by all possible combination, without repetition, of n roots taken at groups of n-1 elements.} can never be divided by 4 (see Appendix B) so, being $\Gamma_1$ and $t_{n-2}$ both even, if $t_{n-1}$ was divisible by 4 then the right side would be a multiple by 4 and therefore also this case is excluded.
\item The term $t_{n-1}$ is the product of the (n-1) roots of the equation (6a) and, as already said, they can be even integers, irrational or complex. In these last two cases they must appear as conjugate pairs\footnote{See Appendix C}.
\end{itemize}
Therefore we have the following two cases:
\begin{description}
\item[a)]	If all the roots $\Gamma_2, \Gamma_3, \ldots \Gamma_{n-1}, \Gamma_n$ are conjugate pairs (irrational or complex) then, also if they fulfill all the relations (8), by Lemma 1 the FLT is proved because there is no any integer root other than the trivial one $\Gamma_1$.
\item[b)]	If some of the $\Gamma_i$  (i=2, 3, \ldots n) were integers then they must be even and at least a pair, therefore  carrying a factor 4 into equation (8c), but this is ruled out by Appendix B.
\end {description}
This exhausts all possible cases, showing that (6a) does not admit integer solutions, other than the trivial one, for any odd integers D and E and for all $n = odd \ primes$ then, by Lemma 1, the Fermat\rq{}s Last Theorem is proved.

\section{Conclusions}
In previous sections we have demonstrated the validity of FLT for n $\geq$ 3 where $A$ e $B$ are both odd (case 1 in Tab.1).

In Appendix A we show that also in cases in which $A$ and $B$ have opposite parity (cases 2 and 3 in Tab.1) the FLT holds.

In conclusion we have proved the validity of Fermat\rq{}s Last Theorem by a procedure, without doubt, Fermat himself could known and then we can infer that this is the {\itshape marvelous proof}, probably been lost, that he claimed to own. The procedure described in this paper does not allow to prove\footnote{See Appendix A} the case n = 4, then we understand why Fermat was worried to demonstrate it in another way.\\

We note that Andrew Wiles proved the FLT only indirectly. In fact Wiles proved the validity of the Taniyama-Shimura conjecture that asserts that every elliptic curve must be related to a modular form. Gerhard Frey had previously devised a mechanism that links the FLT to the elliptic equations and thus indirectly to the Taniyama-Shimura conjecture.

The demonstration of FLT presented in this work, as well as to verify the validity of FLT itself, through the mechanism of Frey, allows us to say that the Taniyama-Shimura conjecture is also verified without the use of the demonstration of Wiles.
\\
%------------------
\appendix

\section{Appendix }
Here, we want analyze the cases 2 and 3 in Tab.1, that is $A $ and $B$ having opposite parity. Of course is enough discuss only the case 2, indeed the case 3 can be reported to case 2 exchanging the variables $A$ and $B$. \\

We start again from relation (1), where now we consider $A=odd$ and $B=even$ and therefore  $C=odd$, then
\begin{eqnarray*} 
A^n+B^n &=& C^n   \  \  \  \  \  \  \  \   or  \hspace{7.3cm} (A1a) \\
A^n-C^n &=& -B^n     \hspace{8.4cm} (A1b) 
\end{eqnarray*}

Let define the two variables\footnote{Similar consideration regarding $D$ and $E$ made into \lq\lq{}First consideration\rq\rq{} paragraph, brings us to conclude that  $GCD (F,G)=1$.}  (similar to $D$ and $E$)
\begin{eqnarray*} 
F=-B-A =odd \ number     \hspace{7cm}(A2a) \\                                                                                                                                                                                                          
G=-B+C=odd \ number    \hspace{7cm}(A2b) 
\end{eqnarray*}
From (A1b), (A2a) and (A2b) we have
$$ 
A^n-C^n = (-B-F)^n  - (B+G)^n = -B^n   \  \  \  \  \  \  \  \  \  \   or  \hspace{4cm} (A3a)
$$
$$
(B+F)^n+(B+G)^n = B^n     \hspace{7.7cm} (A3b) 
$$

\textbf{Note}: the step from (A3a) to (A3b), due to the negative signs inside the first parentheses, can be done only if $n$ is an odd number\footnote{Assuming that the proof given in this paper is actually the {\itshape marvelous proof}  that Fermat claimed to have, probably been lost, then we understand why he worried to demonstrate by other ways the case n = 4.}.\\
\\
Developing the powers of binomials in (A3b), we get:
\begin{eqnarray*}
(B+F)^n &=& \sum_{k=0}^{n} {n \choose k}  F^k B^{n-k}=B^n+\sum_{k=1}^{n} {n \choose k} F^k B^{n-k} \hspace{2.2cm} (A4) \\
(B+G)^n &=& \sum_{k=0}^{n} {n \choose k} G^k B^{n-k}=B^n+\sum_{k=1}^{n} {n \choose k} G^k B^{n-k} \hspace{2.2cm} (A5) \\
\end{eqnarray*}
where, for convenience, we have released the first term under the sign of summation. 
Substituting (A4) and (A5) in (A3b) we obtain the fundamental relationship
\begin{eqnarray*}
P_{(B,n)}&=& B^n+ \sum_{k=1}^{n} {n \choose k} \left ( F^k +G^k\right) B^{n-k}=0  \hspace{4.cm} (A6) \\
\end{eqnarray*}
The (A6) is the expression of a polynomial in the unknown $B$ completely equivalent, except the term $(-1)^k$, to the equation (6) then it leads at an equation similar to (6a), that is:
\begin{eqnarray*} 
P_{(B,n)} & = & B^n+{n \choose 1} \left (F+G\right) B^{n-1}+{n \choose 2} \left ( F^2 +G^2\right) B^{n-2}+... \nonumber \\
                 & + & {n \choose {n-1}} \left ( F^{n-1} +G^{n-1}\right) B+ \left ( F^n +G^n\right)=0  \hspace{2.2cm} (A6a) 
 \end{eqnarray*}                

Now we can borrow all the considerations done on the \lq\lq{}Demonstration of FLT\rq\rq{} paragraph, then showing that the equation (A6a) cannot admit integer roots, other than the trivial one, therefore proving that the Fermat\rq{}s Last Theorem is valid also in the cases 2 and 3 of tab.1.

\section{Appendix }
\newtheorem{Theorem}{Theorem}
\begin{Theorem}
Let $X$ and $Y$ two odd positive integers and n even then the quantity  $X^n+Y^n$ never is divisible by 4.
\end{Theorem}
{\bfseries Proof:}\\
Let n=2, then due to the odd value of $X$, will be either ${X\equiv1\pmod{4}}$  or  ${X\equiv3\pmod{4}}$  , then ${X^2\equiv1\pmod{4}}$   for any odd $X$.
Moreover, ${X^4=X^2X^2\equiv1\pmod{4}}$  so, by induction, ${X^n=X^2X^{n-2} \equiv1\pmod{4}}$    for any even n and odd $X$.\\
To conclude then ${(X^n+Y^n)\equiv2\pmod{4}},$  therefore never divisible by 4.

\section{Appendix }

Actually the equations (8a) and  (8c) pose strong constraints on the values of roots $\Gamma_j$, indeed:

let $\Gamma_{j>1}$ are irrational or integer numbers, then pose $\Gamma_j=\gamma_{j}+\delta_j$ (with $ j>1$) , where $\gamma_{j}$   is the integer part of  $\Gamma_j$  and $0\leq \delta_j < 1$ its decimal irrational part.

So that the sum (8a) is an integer, must be
\[t_1=\sum_{j=2}^{n} \Gamma_{j} = \sum_{j=2}^{n}\left( \gamma_{j} +\delta_j\right)=  \sum_{j=2}^{n} \gamma_{j} +\sum_{j=2}^{n} \delta_j \tag{C1}\label{eq1}\]
where $\displaystyle\sum_{j=2}^{n} \delta_j$  itself must be either integer or null and  $n=odd \ prime$.\\
In similar way from (8c) we have
\[
t_{n-1}=\prod_{j=2}^{n} \Gamma_j= \prod_{j=2}^{n}\left (\gamma_{j} +\delta_j\right)\tag{C2}\label{eq2}
\]
so that both expressions (C1) and  \eqref{eq2} give integer values, needs that the  $\Gamma_j$ have conjugated values at pair\footnote{
%\footnote start
We begin by considering only two terms $\Gamma_j$ and $\Gamma_{j+1}$,  then we must have (from 8a)\\
$S_j=\Gamma_j+\Gamma_{j+1}=
\left(\gamma_{j}+\delta_j\right)+\left(\gamma_{j+1}+\delta_{j+1}\right)=integer$ 
therefore will be $\delta_j+\delta_{j+1}=0$ that is $\delta=\delta_j=-\delta_{j+1}$ and moreover (from 8c)\\
$M_j=\Gamma_j\Gamma_{j+1}=\left(\gamma_{j}+\delta\right)\left(\gamma_{j+1}-\delta\right)=
\gamma_{j}\gamma_{j+1}+\left(\gamma_{j+1}-\gamma_{j}\right)\delta-\delta^2=k$
with $k=integer$  then follows\\
$\delta=\frac{\gamma_{j+1}-\gamma_{j}}{2}\pm\frac{\sqrt{\left(\gamma_{j+1}-\gamma_{j}\right)^2-4\lambda}}{2}$ where $\lambda=k-\gamma_j\gamma_{j+1}$ then getting the positive sign only\\
\\
$\Gamma_j=\left(\gamma_{j}+\delta\right)=\gamma_{j}+\frac{\gamma_{j+1}-\gamma_{j}}{2}+\frac{\sqrt{\left(\gamma_{j+1}-\gamma_{j}\right)^2-4\lambda}}{2}=\alpha_j+\sqrt{\beta_j}$ and\\
\\
$\Gamma_{j+1}=\left(\gamma_{j+1}-\delta\right)=\gamma_{j+1}-\frac{\gamma_{j+1}-\gamma_{j}}{2}-\frac{\sqrt{\left(\gamma_{j+1}-\gamma_{j}\right)^2-4\lambda}}{2}=\alpha_j-\sqrt{\beta_j}$ where\\
\\
$\alpha_j=\frac{\gamma_{j}+\gamma_{j+1}}{2}$ and $\beta_j=\left(\frac{\gamma_{j+1}-\gamma_{j}}{2}\right)^2-\lambda$ \ \ \ therefore\\
\\
$\ S_j=\Gamma_j+\Gamma_{j+1}=2\alpha_j=\gamma_j+\gamma_{j+1}$ \\
$M_j=\Gamma_j\Gamma_{j+1}={\alpha_j}^2-\beta_j=\gamma_j\gamma_{j+1}+\lambda$ \\

Taking in account more terms $\Gamma_i$ (i=4, 6,...n) we obtain similar results where the $\alpha_i$ and $\beta_i$ will be function of the corresponding $\gamma_i$, always taken in pair.

}, i.e.:
%\footnote end
\\
$\Gamma_j=\alpha_j+\sqrt{\beta_j} $\\
$\Gamma_{j+1}=\alpha_j-\sqrt{\beta_j}$ \\

Where $\alpha_j=\frac{\gamma_{j}+\gamma_{j+1}}{2}$ and  $\beta_j=\left[\frac{\gamma_{j+1}-\gamma_{j}}{2}\right]^2-\lambda$  with $\lambda$ a suitable integer and $(j=2,4,...n)$.\\

\begin{flushleft}

\end{flushleft}

\end{document}